\documentclass[11pt]{amsart}

\title{A Polynomial Time Algorithm for the Braid Double Shielded Public Key Cryptosystems  }

\author{Vitali\u\i\ Roman'kov}
\address{Institute of Mathematics and Information Technologies\\Omsk State Dostoevskii University}
\curraddr{}
\email{romankov48@mail.ru}

\usepackage{amsmath,amsthm,amsfonts,amssymb}
\usepackage{graphicx}
\usepackage{hyperref}
\usepackage[usenames]{color}

\theoremstyle{definition}

\newcounter{comcount}

\def\F{{\mathbb{F}}}

\date{}

\begin{document}

\maketitle

\begin{abstract}
We propose new provable practical deterministic polynomial time algorithm for the braid Wang, Xu, Li, Lin and Wang {\it Double shielded public key cryptosystems} \cite{WXLLW}. We show that a linear decomposition attack based on the decomposition method introduced by the author in  monography \cite{R1} and  paper \cite{R2} works for the image of braids under the Lawrence-Krammer representation by finding the exchanging keys in the both two main protocols in \cite{WXLLW}. 
\end{abstract}

\section{Introduction}
\label{se:intro}

In this paper we discuss, following \cite{R1} and \cite{R2} a new practical  attack on  the two main protocols proposed in \cite{WXLLW}.  This attack introduced by the author in \cite{R1} and \cite{R2} works when the platform  groups  are linear. We show that in this case, contrary to the common opinion (and some explicitly stated security assumptions), one does not need to solve the underlying algorithmic problems to break the scheme, i.e.,  there is another  algorithm that recovers the private keys without solving the principal algorithmic problem on which the  security assumptions are based. This changes completely our understanding of security of these scheme. The efficacy of the attack depends on the platform group, so it requires a  specific analysis in each particular case. In general one can only state that the attack is in polynomial time in the size of the data, when the platform and related groups are given together with their linear representations. In many other cases we can effectively use known linear presentations of the groups under consideration. The braid groups are among them in view of  the Lawrence-Krammer representation.  

The monography \cite{R1} solves uniformly  protocols based on the conjugacy search problem (Ko, Lee et. al. \cite{KL}, Wang, Cao et. al \cite{WC}), protocols based on the decomposition and factorization problems (Stickel \cite{S}, Alvares, Martinez et. al. \cite{AM}, Shpilrain, Ushakov \cite{SU}, Romanczuk, Ustimenko \cite{RU}), protocols based on actions by automorphisms (Mahalanobis \cite{M}, Habeeb, Kahrobaei et. al. \cite{HK}, Markov, Mikhalev et. al. \cite{MM}), and a number of other protocols. 

In this paper we apply our method to the double shielded key exchange protocols 1 and 2 proposed in \cite {WXLLW}.

\section{Constructing a basis}

This section is a part of (still unpublished) paper of the author and A.G. Miasnikov. 

Let $V$ be a finite dimensional vector space over a field $\mathbb{F}$  with basis $\mathcal{B} = \{v_1, \ldots,v_r\}$.  Let $End(V)$ be  the  semigroup of endomorphisms of $V.$  We assume that  elements $v \in V$ are given as vectors relative to $\mathcal{B}$, and endomorphisms $a \in End(V)$ are given by their matrices relative to $\mathcal{B}$. Let $<W>$ denotes submonoid generated by $W.$ 
For an endomorphism $a \in End(V)$ and an  element $v \in V$ we denote by  $v^a$ the image of $v$ under $a.$  Also, for any subsets $W \subseteq V$ and $A \subseteq End(V)$ we put $W^A = \{w^a | w \in W, a \in A\}$, and  denote by Sp$(W)$ the subspace of $V$ generated by $W$. We assume that elements of the field $\F$ are given in some constructive form and the "size" of the form is defined. Furthermore, we assume that the basic field operations in $\F$ are efficient, in particular they can be performed in polynomial time in the size of the elements. In all the particular protocols considered in this paper  the field $\F$ satisfies all  these conditions. 

There is an algorithm that for given finite subsets $W \subseteq V$ and $U  \subseteq End(V)$ finds  a basis of the subspace Sp$(W^{\langle U\rangle})$   in the form $w_1^{a_1}, \ldots, w_t^{a_t}$, where $w_i  \in W$ and $a_i$ is a  product  of elements  from $U$.  Furthermore, the number of field operations used by the algorithm is polynomial in $r = \dim_\F V$ and the cardinalities of  $W$ and $U$.

Using Gauss elimination one can effectively find  a maximal linearly independent subset $L_0$ of $W$. Notice that Sp$(L_0^{\langle U\rangle}) =$  Sp$(W^{\langle U\rangle})$. Adding to the set $L_0$  one by one elements $v^a$, where $v \in L_0, a \in U$ and checking every time  linear independence of the extended set, one can effectively construct a maximal linearly independent subset $L_1$ of the set $L_0 \cup L_0^U$ which extends the set $L_0$. Notice that Sp$(L_0^{\langle U\rangle})$ = Sp$(L_1^{\langle U\rangle})$ and the elements in $L_1$ are of the form $w^a$, where $w \in W$ and $a = 1$ or $a \in U$. It follows that if $L_0 = L_1$ then $L_0$ is a basis of Sp$(W^{\langle U\rangle})$. If $L_0 \neq L_1$ then we repeat the procedure for $L_1$ and find a maximal linearly independent subset $L_2$ of $L_1 \cup L_1^U$ extending $L_1$. Keep going one constructs a sequence of strictly increasing subspaces $L_0 < L_1 < \ldots < L_i$ of $V$. Since the dimension $r$ of $V$is finite  the sequence stabilizes for some $i \leq r$. In this case  $L_i$ is a basis of  Sp$(W^{\langle U\rangle})$ and its elements are in the required form. 

To estimate the upper bound of the number of the field operations used by the algorithm, observe first that the number of the field operations in Gauss elimination performed on  a matrix of size $n\times r$ is $O(n^2r)$. Hence it requires at most $O(n^2r)$ steps to construct $L_0$ from $W$, where $n = |W|$ is the number of elements in $W$. Notice that $|L_j| \leq r$ for every $j$. So to find $L+{j+1}$ it suffices to perform Gauss elimination on the matrix corresponding to $L_j \cup L_j^U$ which has size at most $r +r|U|$.Thus the upper estimate on this number is $O(r^3|U|^2)$. Since there are at most $r$ iterations of this procedure one has the total estimate as $O(r^3|U|^2 +r|W|^2)$.

In this paper $V$ is underlying linear space of a matrix algebra Mat$_t(\mathbb{F})$ of all matrices of size $t \times t$ over $\mathbb{F}.$ Let $G$ be a subgroup of the multiplicative group of Mat$_t(\mathbb{F})$, and $A$ and $B$ are two subgroups of $G.$ Every pair of elements $a \in A$ and $b \in B$ define an automorphism $\varphi_{a, b}$ of $V$ such that for every $v \in V$ one has  $v^{\varphi_{a,b}} = avb.$ Let $U$ be submonoid generated by all such automorphisms. Thus for every subset $W \subseteq V$ one can effectively construct a basis of subspace $W^{U}.$

\section{The double shielded key exchange protocol  1 \cite{WXLLW}}.

In view of the Lawrence-Krammer representation of the braid group $B_n$ we can assume that the group $G$ below is given as a linear group over a field $\mathbb{F}$. So, we assume that $G$ is a part of a finite dimensional vector space $V.$ 

At first we describe  the protocol 1.

\begin{itemize}
\item Alice and Bob agree on a nonabelian group $G,$ and randomly chosen element $h \in G$ and two subgroups $A, B$ of $G$, such that $ab=ba$ for any $a \in A$ and any $b \in B.$ We assume that $A$ and $B$ are finitely generated and are given by the fixed generating sets $\{a_1, ..., a_n\}$ and $\{b_1, ..., b_m\},$ respectively. 
\item Alice chooses four elements $c_1, c_2, d_1, d_2  \in A,$ computes $x = d_1c_1hc_2d_2,$ and then sends $x$ to Bob.
\item Bob chooses six elements $f_1, f_2, g_1, g_2, g_3, g_4 \in B$, computes $y = g_1f_1hf_2g_2$ and $w = g_3f_1xf_2g_4$, and then sends $(y, w)$ to Alice.
\item Alice chooses two elements $d_3, d_4 \in A$, computes $z =  d_3c_1yc_2d_4$ and $u = d_1^{-1}wd_2^{-1}$, and then sends $(z, u)$ to Bob. 
\item Bob sends $v = g_1^{-1}zg_2^{-1}$ to Alice. 
\item Alice computes $K_A = d_3^{-1}vd_4^{-1}$.
\item Bob computes $K_B = g_3^{-1}ug_4^{-1} = c_1f_1hf_2c_2$ which is equal to $K_A$ and then $K = K_A = K_B$ is Alice and Bob's common secret key. 
\end{itemize}

Now we show how the common secret key can be computed. 

Let  $BzB$ be  subspace of $V$ generated by all elements of the form $fzg$ where $f, g \in B.$ We can construct a basis $\{e_izl_i (e_i, l_i \in B, i = 1, ..., r)\}$ of $BzB$ in a polynomial time as it is  explained in the previous section. Since $v \in BzB,$ we can effectively write it in the form 

\begin{equation}
\label{eq:1}
v = \sum_{i=1}^r\alpha_i e_izl_i,
\end{equation}

\noindent
where $\alpha_i \in \mathbb{F}, i=1, ..., r.$

In a similar way we construct  bases $\{e_j'hl_j' (e_j', l_j' \in B, j = 1, ..., s)\}$ of  $BhB$, and  $\{e_k''wl_k''
(e_k'', l_k'' \in B, k = 1, ..., q)\}$ of $BwB.$ Then we get presentations

\begin{equation}
\label{eq:2}
y = \sum_{j=1}^s\beta_j e_j'hl_j',
\end{equation}

\noindent
where $\beta_j \in \mathbb{F}, j=1, ..., s,$

\noindent and
\begin{equation}
\label{eq:3}
x = \sum_{k=1}^q\gamma_k e_k''wl_k'',
\end{equation}

\noindent
where $\gamma_k \in \mathbb{F},  k=1, ..., q.$

Now we swap $w$ by $u$ in the right hand side of (\ref{eq:3}). By direct computation we obtain 

\begin{equation}
\label{eq:4}
\sum_{k=1}^q\gamma_k e_k''ul_k'' = \sum_{k=1}^q\gamma_k e_k''d_1^{-1}wd_2^{-1}l_k'' = d_1^{-1}(\sum_{k=1}^q\gamma_k e_k''wl_k'')d_2^{-1} =d_1^{-1}xd_2^{-1} = c_1hc_2.
\end{equation}

Then we swap $h$ by $c_1hc_2$ in the right hand side of (\ref{eq:2}). 

\begin{equation}
\label{eq:5}
 \sum_{j=1}^s\beta_j f_j''c_1hc_2g_j'' = 
c_1(\sum_{j=1}^s\beta_j e_j'hl_j')c_2=c_1yc_2=
c_1g_1f_1hf_2g_2c_2. 
\end{equation}

At last we swap $z$ by $c_1g_1f_1hf_2g_2c_2$ 
in the right hand side of (\ref{eq:1}) and get 

\begin{equation}
\label{eq:6}
\sum_{i=1}^r\alpha_i e_ic_1g_1f_1hf_2g_2c_2l_i = d_3^{-1}(\sum_{i=1}^r\alpha_i e_izl_i)d_4^{-1} = c_1f_1hf_2c_2 = K.
\end{equation}

\section{The double shielded key exchange protocol  2 \cite{WXLLW}}.

As before we  assume that the group $G$ below is given as a linear group over a field $\mathbb{F}$. So, we assume that $G$ is a part of a finite dimensional vector space $V.$ 

Firstly we describe  the protocol 2.

\begin{itemize}
\item Alice and Bob agree on a nonabelian group $G,$ and randomly chosen element $h \in G$ and two subgroups $A, B$ of $G$, such that $ab=ba$ for any $a \in A$ and any $b \in B.$ We assume that $A$ and $B$ are finitely generated and are given by the fixed generating sets $\{a_1, ..., a_n\}$ and $\{b_1, ..., b_m\},$ respectively. 
\item Alice chooses four elements $c_1, d_1 \in A$ and $f_2, g_2 \in B,$  computes $x = d_1f_1hf_2g_2,$ and then sends $x$ to Bob.
\item Bob chooses six elements $c_2, d_2, d_3 \in A$ and $ f_1, g_1, g_3 \in B$, computes $y = g_1f_1hc_2d_2$ and $w = g_3f_1xc_2d_3$, and then sends $(y, w)$ to Alice.
\item Alice chooses two elements $ d_4 \in A$ and $g_4 \in B$,  computes $z =  d_4c_1yf_2g_4$ and $u = d_1^{-1}wg_2^{-1}$, and then sends $(z, u)$ to Bob. 
\item Bob sends $v = g_1^{-1}zd_2^{-1}$ to Alice. 
\item Alice computes $K_A = d_4^{-1}vg_4^{-1}$.
\item Bob computes $K_B = g_3^{-1}ud_3^{-1} = c_1f_1hf_2c_2$ which is equal to $K_A$ and then $K = K_A = K_B$ is Alice and Bob's common secret key. 
\end{itemize}

Now we show how the common secret key can be computed. 

Let  $BzA$ be  subspace of $V$ generated by all elements of the form $fzd$ where $f \in B, d \in A.$ We can construct a basis $\{e_izl_i (e_i \in B, l_i \in A, i = 1, ..., r)\}$ of $BzA$ in a polynomial time as it is  explained in the previous section. Since $v \in BzA,$ we can effectively write 

\begin{equation}
\label{eq:7}
v = \sum_{i=1}^r\alpha_i e_izl_i,
\end{equation}

\noindent
where $\alpha_i \in \mathbb{F}, i=1, ..., r.$

In a similar way we construct  bases $\{e_j'hl_j' (e_j' \in B, l_j' \in A, j = 1, ..., s)\}$ of  $BhA$, and  $\{e_k''wl_k''
(e_k'' \in B, l_k'' \in A, k = 1, ..., q)\}$ of $BwA.$ Then we write

\begin{equation}
\label{eq:8}
y = \sum_{j=1}^s\beta_j e_j'hl_j',
\end{equation}

\noindent
where $\beta_j \in \mathbb{F}, j=1, ..., s,$

\noindent and
\begin{equation}
\label{eq:9}
x = \sum_{k=1}^q\gamma_k e_k''wl_k'',
\end{equation}

\noindent
where $\gamma_k \in \mathbb{F},  k=1, ..., q.$

Now we swap $w$ by $u$ in the right hand side of (\ref{eq:3}). By direct computation we obtain 

\begin{equation}
\label{eq:10}
\sum_{k=1}^q\gamma_k e_k''ul_k'' = \sum_{k=1}^q\gamma_k e_k''d_1^{-1}wg_2^{-1}l_k'' = d_1^{-1}(\sum_{k=1}^q\gamma_k e_k''wl_k'')g_2^{-1} =d_1^{-1}xg_2^{-1} = c_1hf_2.
\end{equation}

Then we swap $h$ by $c_1hf_2$ in the right hand side of (\ref{eq:2}). 

\begin{equation}
\label{eq:11}
 \sum_{j=1}^s\beta_j e_j'c_1hf_2l_j' = 
c_1(\sum_{j=1}^s\beta_j e_j'hl_j')f_2=c_1yf_2=
c_1g_1f_1hc_2d_2f_2. 
\end{equation}

At last we swap $z$ by $c_1g_1f_1hc_2d_2f_2$ 
in the right hand side of (\ref{eq:1}) and get 

\begin{equation}
\label{eq:12}
\sum_{i=1}^r\alpha_i e_ic_1g_1f_1hc_2d_2f_2l_i = d_4^{-1}(\sum_{i=1}^r\alpha_i e_izl_i)g_4^{-1} =  d_4^{-1}vg_4^{-1} = c_1f_1hc_2f_2 = K.
\end{equation}

Two other, the shielded public key encryption protocol and the shield digital signature protocol in \cite{WXLLW} completely based on the protocols 1 and 2.

\section{The Lawrence-Krammer Representation}

Let $B_n$ denote the Artin braid group on $n$ strings, $n \in \mathbb{N}.$ R. Lawrence described in 1990 a family of so called Lawrence representations of $B_n.$ Around 2001 S. Bigelow \cite{Big} and D. Krammer \cite{Kr} independently proved that all braid groups $B_n$ are linear. Their work used the  Lawrence-Krammer representations $\rho_n : B_n \rightarrow GL_{n(n-1)/2}(\mathbb{Z}[t^{\pm 1}, s^{\pm 1}])$ that   has been proved faithful for every $n\in \mathbb{N}.$  One can effectively find the image $\rho_n (g)$ for every element $g \in B_n.$ 
Moreover, there exists an effective procedure to recover a braid $g \in B_n$ from its image $\rho_n (g).$ It was shown in \cite{CJ} that it can be done in $O(2m^3log d_t)$ multiplications of entries in $\rho_n (g).$ Here $m = n(n-1)/2$ and $d_t$ is a parameter  that can be effectively computed by $\rho_n (g).$  See \cite{CJ} for details.

\section{Conclusion}

In this paper we proposed a polynomial time deterministic algorithm to recover secret keys established by the protocols 1 and 2 in \cite{WXLLW}. We assumed that the group $G$ in this protocols is linear.

The authors of \cite{WXLLW} suggested that the infinite nonabelian groups   $B_n$ with $n \geq 12$ can be taken as the platform groups for  the protocols 1 and 2 in the above  section 2. By the Lawrence-Krammer representations the groups $B_n$ are linear. Moreover, this representations are effective computable and invertible.  Unfortunately, in this setting the proposed protocols are not secure. Our  cryptanalysis in the above section 2 shows that the linear decomposition attack  works effectively in this case.

\end{document}